\documentclass[12pt]{article}
\usepackage{latexsym}
\usepackage{amssymb}

\newtheorem{teo}{THEOREM}[section]
\newtheorem{prop}[teo]{PROPOSITION}
\newtheorem{lem}[teo] {LEMMA}

\newtheorem{obser}[teo]{Remark}
\newtheorem{defi}[teo]{DEFINITION}
\newtheorem{coro}[teo]{COROLLARY}
\newenvironment{note}{\noindent Note: \rm}

\newenvironment{Demo}{\noindent \sf Proof. \rm}

\def\Hom{\mathop{\rm Hom}\nolimits}

\def\qed{\hfill \mbox{$\square$}}

\def\Ext{\mathop{\rm Ext}\nolimits}

\def\Coker{\mathop{\rm Coker}\nolimits}

\def\dim{\mathop{\rm dim_k}\nolimits}

\def\centro{\mathop{\rm Center}\nolimits}

\def\ind{\mathop{\rm ind}\nolimits}

\def\ra{\rightarrow}

\begin{document}

\date{}

\title{Hochschild Cohomology of skew group rings and invariants
\thanks{ This work was done during the visit of the second author at the Universidade de S\~ao Paulo on December 2000 and it is part of the work supported by a interchange grant from CNPq (Brasil) and CONACyT(Mexico). The 
first and the second authors thank respectively CNPq and CONACyT.}}

\maketitle

\centerline{\author{E. N. Marcos, R. Mart\'\i nez-Villa and Ma. I. R. Martins}}

\vspace*{.5cm}

\begin{center}
{{\em{\rm In memory of Sheila Brenner}}}
\end{center}

\begin{abstract}
Let $A$ be a $k$-algebra and $G$ be a group acting on $A$. We show that  $G$ also acts on  the Hochschild cohomology algebra 
$HH^{\bullet}(A)$ and that  there is a monomorphism of rings $HH^{\bullet}(A)^G \hookrightarrow HH^{\bullet}(A[G])$. That allows us to show the existence of a monomorphism from  $HH^{\bullet}(\widetilde{A})^G$ into $HH^{\bullet}(A)$, where $\widetilde{A}$ is a Galois covering with group $G$.
\end{abstract}



\section {\sf Introduction }
 The Hochschild cohomology groups were introduced by Hochschild fifty years ago, but they have been investigated lately under different aspects
by many authors. 

In this work our interest is to study the Hochschild cohomology of  skew group rings and of certain Koszul algebras. Our main purpose is comparing the Hochschild cohomology algebra $HH^{\bullet}(A)=
\oplus_{i\geq 0}\Ext^i_{A^e}(A, A)$ of a $k$-algebra $A$ with the Hochschild cohomology algebra $HH^{\bullet}(A[G])$ of the skew group algebra $A[G]$, with $G$ being a finite group acting on $A$. We obtain the following relation between these cohomology algebras.

{\bf Theorem  2.9. } {\it  Let $A$ be a $k$-algebra and $G$ be a finite group acting on $A$. Then $G$ acts on the Hochschild cohomology $k$-algebra $HH^{\bullet}(A)$, and there is a ring monomorphism: $HH^{\bullet}(A)^G\hookrightarrow HH^{\bullet}(A[G]).$}

 For  finite groups, it is known that there is a strong connection between skew group rings and Galois coverings, and between skew group rings  and smash products of graded algebras (see \cite{CM}). These facts  and  the existence of the monomorphism  in Theorem above lead us to investigate a possible relation between the Hochschild cohomology algebras of $G$-graded algebras and their covering algebras defined by $G$. In this direction,  let  $A$ be a $G$-graded algebra, with $G$ being a finite group.  We consider the covering algebra of $A$ defined by $G$ and we relate it with the smash product $A\# kG^*$. Then, as a consequence of Theorem \ref{main} and using Theorem of duality coactions, we obtain  that $G$ is a group of automorphism of $\widetilde{A}$ and the existence of a ring monomorhism from
$HH^{\bullet}(\widetilde{A})^G$ into $HH^{\bullet}(A)$.

While we were reviewing the final version of this article we received from C. Cibils and M.J. Redondo a pre-print  entiltled  {\it Cartan-Leray spectral sequence for Galois coverings of categories}, \cite{cr}. In this pre-print they gave a spectral sequence involving the Hochschild cohomology  of  Galois coverings and they show that the monomorphism obtained by us  is an isomorphism in case the characteristic of the field is zero.

 In this work we also study the Hochschild cohomology groups of Koszul algebras of finite global dimension. With this we reach the  Hochschild cohomology groups  of $\mathbb{C}$-preprojective algebras associated to Euclidean diagrams and of Auslander algebra of standard algebras, since these algebras are examples of Koszul algebras of global dimension two.  We obtain a lower bound for dimension of $HH^n(A)$, for $A$ be a Koszul algebra of global dimension $n$. So as consequence  we get that the second Hochschild cohomology group of  $\mathbb{C}$-preprojective algebra of Euclidean type  does not vanish; and  it is also true for Auslander algebra of standard algebras having  non projective indecomposable modules isomorphic to their own Auslander-Reiten translate.  So,  in both these cases the algebras are not rigid (see \cite{ge})

We now describe the contents of each section  in the paper. In  section 2, after recalling some notions and known facts needed along the work, we state and prove the main result of the section - Theorem \ref{main}.    
 This theorem states that  if a group $G$ is a group acting  on an algebra $A$,  then $G$ also acts on the Hochschild cohomology algebra of $A$ and there is a ring monomorphism between the fixed points  of the Hochschild cohomology algebra of $A$ and the Hochschild cohomology algebra  of the skew group ring $A[G]$.

In section \ref{covering} we define the  covering algebra $\widetilde{A}$ associated to a $G$-graded algebra $A$, where $G$ is finite group. We also recall the notion of smash product, and we show that this product  is isomorphic to that covering algebra. This  isomorphism together with   Theorem \ref{main} and duality coactions gives us a  similar relationship between the invariants of  Hochschild cohomology ring $HH^{\bullet}(\widetilde{A})$ and  the Hochschild cohomology ring of $A$.

In section \ref{Koszul} we deal with quadratic algebras. We construct the Koszul complex for  quotient of path algebras by quadratic ideals through a similar procedure  used by Berger in \cite{Be}. The Koszul complex (named bimodule Koszul complex by him) was constructed in \cite{Be} for quotient of free associative algebras $A$ and it is  minimal graded resolution of $A$ as $A-A$ bimodule, in case $A$ is a generalized Koszul algebra (also called $d$-Koszul algebras). But it can be also constructed for algebras  which are quotient of quiver algebras by ideals generated by elements of degree $d\geq 2$ (see \cite{gmmp}). Our construction here follows closely the one in \cite{gmmp} for Koszul algebras (that is, 2-Koszul). It enables us  to obtain a lower bound for the dimension of the
$n$-Hochschild cohomology group of Koszul algebras of global dimensional $n$, and in consequence the property mentioned above for preprojective algebras of Euclidean-type and Auslander algebra of an algebra $A$ standard.

\section {\sf Hochschild cohomology rings and invariants }\label{basics}

Given a ring $A$ we denote by $A^{op}$ its opposite ring.
For $a\in A$ we denote by $a^o \in A^{op}$ the corresponding  element  in $A^{op}$. In case  that $A$ is an algebra over a field $k$ we will denote by $A^e$ its enveloping algebra $A\otimes_k A^{op}$. Moreover, if $A$ and $B$ are algebras over $k$, the algebra  $A\otimes_{_k} B$ will be denoted simply by $A\otimes B$. Sometimes, by
 simplicity,  we will not explicit the ground ring of tensor product when it is clear in the context. 

We also remark that the category of left modules over the algebra $A^e$ is canonically isomorphic to the category of $A-A$ bimodules.  So we use this isomorphism as identification.

Now we recall some definitions and basic facts.

\begin{defi}
Let $A$ be a ring and $G$ a group. We say that $G$ acts on $A$   if there is a group homomorphism between $G$ and the group Aut$(A)$ of ring automorphism of $A$. If this group homomorphism is injective, we say that $G$ acts faithfully on $A$ or that $G$ is a group of automorphism of $A$.

\end{defi}

We remark that if $G$ acts on $A$,  then $G$ naturally  also acts on the opposite ring $A^{op}$.  In  case that $A$ is a $k$-algebra, we will assume that the group Aut$(A)$ is the group of  automorphisms of $k$-algebras. Moreover, if $G$ and $H$ are groups acting on the $k$-algebras $A$ and $B$, respectively, then the group $G\times H$ acts on $A\otimes_k B$, and consequently 
 $G\times G$ acts on $A \otimes A^{op}$.

 Now we are going to recall the definition of {\it skew group algebra}.

\begin{defi}

Let $A$ be a ring and $G$ a finite group acting on $A$. The elements of the skew group ring $A[G]$ are the same as those of the corresponding group ring. Addition is as usual coordinate-wise, and multiplication is extended by bilinearity  from the formula $(ag)(bh)= a g(b)gh$, for $a$ and $b$ in $A$ and $g$ and $h$ in $G$.

\end{defi}

The following statements will be useful later and their verification are routine. 

\begin{prop}\label{isotimes} Let $A$ be a $k$-algebra and $G$ be a finite group acting on $A$. Then,
\begin{enumerate}
\item[i)] The rings $(A[G])^{op}$ and  $A^{op} [G]$ are isomorphic, via the map \linebreak $\theta( (ag)^o ) = (g^{-1}(a))^o g^{-1}$.
\item[ii)]  $A[G]\otimes_k (A[G])^{op}$ is isomorphic to $(A\otimes_k A^{op})[G\times G]$, via the map \linebreak$\psi(ag\otimes (bh)^o)= a\otimes (h^{-1}(b))^o (g, h^{-1})$.
\end{enumerate}
\end{prop}

Now we are going to describe certain approach on the category of the left $A$-modules and the category of the left $A[G]$-modules, where $G$ is a finite group acting on the algebra $A$. This approach will be very useful for the next sections.
 
Denoting by Mod-$A$ the category of the left $A$-modules, for each $g\in G$, we can associate a functor, denoted by $^g ( )$, on Mod-$A$. This functor associates to each $M$ in Mod-$A$ the module $^gM$ defined as follows: $^gM= M$ as an abelian group (or k-vector space, in case $A$ is a $k$-algebra) and for $a\in A$ and $m\in M$, $a\cdot_{_g} m:= g(a)m$. On the morphism, $^g ( )$ is defined as the identity.  Observe that the functor $^g ( )$ is clearly an exact functor and is an automorphism of Mod-$A$. We also observe that it is possible to define in a similar faction  an automorphism $( )^g$ on the category $A$-Mod of the right $A$-modules. Furthermore, analogously  we could be consider a functor $ ^g ( )$ on the category of $A- A$ bimodules, by considering $^g M = M$ as an $A- A$ bimodule where it has on the left the structure as above and on the right the original  structure of $M_A$ (also in a similar way we could have $M^g$ as a $A- A$ bimodule).

The following facts can be verify easily.

\begin{prop}\label{Misogm0} Let $A$ be a ring and $G$ a group acting on $A$. If $g$ and $h$ are in $G$ and $M$ is a left $A$-module, then:

\begin{enumerate}

\item[i)] $^g(^hM)=\, ^{hg}\!M$;

\item[ii)]$(M^g)^h = M^{gh}$;

\item[iii)]$ Ag \cong A^g \cong \,\,  ^{g^{^{\!\!\!-1}}}\!\!\!\!A$ (as $A-A$ bimodules).

\end{enumerate}

\end{prop}

Now we recall some basic facts related to $A[G]$-modules (see \cite{M})

\begin{prop}\label{Misogm}

Let $A$ be a ring and $G$ a finite group acting on $A$. If $M$ is in Mod-$A[G]$, then  the map $^g \Psi: M \to ^g\!\!\!M$ given by
$^g \Psi(m) = g(m)$ defines an  isomorphism of $A$-modules.

\end{prop}

\begin{Demo} Clearly $^g \Psi$ preserves the sum. Let $a\in A$ and $m\in M$. So, we have $^g \Psi(am) = g(am)= g(a)g(m) = a\cdot\!\!\!_{_g}\,g(m)= a\cdot\!\!_{_g} \, ^{g}\Psi(m)$, what  shows that $^g\Psi$ is a homomorphism of $A$-modules. A similar verification shows that the map $ m \to g^{^{-1}}(m)$ is the $A$-morphism inverse of $^g \Psi$.
\qed

\end{Demo}

 \begin{prop}\label{gmod}(Lemma 4 in \cite{M})
Let $A$ be a $k$-algebra and $G$ a finite group acting on $A$. Let $M$ and $N$ be $A[G]$-modules. Then the following statements hold:

\begin{description}

\item[i)] The abelian group $\Hom_A(M, N)$ is a $kG$-module, with the action  $(g*f)(m) = g(f(g^{-1}(m))$. Denote by $\Hom_A(M, N)^G$ the set of fixed points, then\linebreak 
 $\Hom_A(M, N)^G =\Hom_{A[G]}(M, N)$;

\item[ii)]For all $i\geq 1$, there is a natural action of $G$ on $\Ext^i_A(M, N)$ and it verifies that $\Ext^i_{A[G]} (M, N) = \Ext^i_A (M, N) ^G$. If $g\in G$, $\xi\in \Ext^i_A(M,N)$ and $\eta\in \Ext^j_A(M, N)$, then $g(\xi\eta)=g (\xi) g(\eta)$.

\end{description}

\end{prop}

\begin{Demo}

i)  It is easy and well known (see for instance \cite{M}).

ii) Let $g\in G$. Since $^{g^{-1}}()$ is an exact functor, then a given element $\xi\in\Ext^i_A(M, N)$ is taken to an element $^{g^{-1}}\!\xi\in \Ext^i_A(^{g^{^{-1}}}\!\!M,\,  ^{g^{^{-1}}}\!\!N)$. But,  since $M$ and $N$ are $A[G]$-modules, by using 
  the isomorphism $^{g^{-1}}\Psi$ and its inverse (see Proposition \ref{Misogm})  we get an  exact sequence, denoted by $g(\xi)$, which is an element in $\Ext^i_A(M , N)$. We note that if two exact sequences are representives of the same element in $\Ext^i_A(M, N)$, then their correspondents under $g(-)$ have the same property. Hence  it  indicates how to define the action. Proposition \ref{Misogm0} guarantees that it really defines an  action of $G$  on $Ext^i_A(M,N)$.

The rest of the proof follows from the fact that the functor $^{g^{-1}}()$ also preserves the $A$-projective modules  and we leave the details to the reader.
\qed
\end{Demo}

\begin{note}
We note that statement ii) in the last proposition could also be proved by remarking that $A[G]$ is a projective $A$-module, and applying the functor $\Hom_A(-, N)$ to the projective resolution of $M$ as an $A[G]$-module, and observing that the homology obtained at each step is a $kG$-module. 

\end{note}

Now we are going to recall the definition of Hochschild cohomology. We present it closely to the original approach given by Hochschild  and it can be found, for instance, in \cite{ibra, caei, H,  W}.

\begin{defi}

Let $A$ be an algebra over a commutative ring $R$, and $_AX_A$ be an $A$-bimodule. The $i^{th}$-Hochschild cohomology group of $A$ with coefficients in $X$, denoted by $HH^i(A, X)$, is the $i^{th}$ cohomology group of the following complex: 

$$0\rightarrow X \stackrel{d_0}{\rightarrow} \Hom_R(A,X) \stackrel{d_1}{\rightarrow}\Hom_R(A^{\otimes 2}, X) \cdots
 \Hom_R(A^{\otimes n},X) \stackrel{d_n}{\rightarrow} \Hom_R(A^{\otimes n+1}, X) \cdots $$ where, for
$n\geq 1$ and $f\in \Hom_R(A^{\otimes n},X)$, $d_n$ is defined by

\begin{eqnarray*}
d_n(f)(x_1 \otimes x_2\dots \otimes x_{n+1})&=& x_1 f(x_2\otimes \dots\otimes x_{n+1})\\ &+& \sum_{i=1}^n (-1)^i f(x_1
\otimes \dots \otimes x_ix_{i+1} \otimes \dots \otimes x_{n+1})\\ &+& (-1)^{n+1} f(x_1 \otimes \dots \otimes
x_{n})x_{n+1}
\end{eqnarray*}

\noindent and for $x\in X$ and $a\in A$, $(d_0 x)(a) =  ax - xa .$

\end{defi}

In  case that $R=k$ is  a field, a different way of approaching to Hochschild  cohomology groups is to consider the  enveloping algebra $A^e=A\otimes_kA^{op}$. In this case $HH^i(A, X) = Ext_{A^e}^i(A, X)$, for all $i\geq 0$. But our particular interest is the example $X = A$, whose $HH^i(A, A)$ is simply denoted by $HH^i(A)$, for $i\geq 0$. These groups are used to define the {\it Hochschild cohomology algebra} $HH^{\bullet}(A)=\oplus_{i\geq 0}HH^i(A)= \oplus_{i\geq 0}\Ext^i_{A^e}(A, A)$ with the multiplication induced by the Yoneda product. In this way $HH^{\bullet}(A)$ is a $\mathbb{Z}$-graded algebra (see for instance \cite{ibra,  W}).

The facts  which we state next point up how important they are to stablish the main result  of this section:  to relate the Hochschild cohomology algebras $HH^{\bullet}(A)$ and $HH^{\bullet}(A[G])$.

Let $A$ be a $k$-algebra and $G$ be a finite group acting on $A$. We have seen, in  Proposition \ref{isotimes}, the enveloping algebra $A[G]^e$ of the skew group algebra $A[G]$ is isomorphic to the $k$-algebra $A^e[G\times G]$. We shall  use this isomorphism together with Proposition \ref{gmod} for describing  an action of $G\times G$ on the Hochschild cohomology algebra $HH^{\bullet}(A[G])$, which respects its    $\mathbb{Z}$-grading (meaning  the action takes an element in $HH^i(A[G])$ to an element in $HH^i(A[G])$. 

With this in mind, we recall that  $HH^{\bullet}(A)=\oplus_{i\geq 0}\Ext^i_{A^e}(A, A)$ and 
 $HH^{\bullet}(A[G])=\oplus_{i\geq 0}\Ext^i_{A[G]^e}(A[G], A[G])$. We also recall that $A[G]$ is a $A^e[G\times G]$-module, with the following `` multiplications '':
$(x\otimes y^o)(\sum_{g\in G} a_g g)= \sum_{g\in G} xa_gg(y) g$, and $(\sigma, \tau)(\sum_{g\in G} a_g g)= 
\sum_{g\in G}\sigma(a_g)\sigma g\tau^{-1}$, for $x$, $y$, $a_g$ in $A$, and $\sigma, \tau$ in $G$.
Then from Proposition \ref{gmod}, for each $i\geq 0$, follows that $G\times G$ acts on $\Ext_{A^e}^i(A[G], A[G])$ and $\Ext^i_{A[G]^e}(A[G], A[G])\cong(\Ext^i_{A^e}(A[G], A[G])^{^{G\times G}}$. So, from the action of $G\times G$ on the grading  we   get the one wanted  on $HH^{\bullet}(A[G])$.

For obtaining the relationship between the  Hochschild cohomology algebras of $A$ and of $A[G]$, now we describe a bit more in details the action considered above (the one considered  in Proposition 
\ref{gmod}) just for the cases $i=0$ and $i= 1$, since for $i\geq 2$ the procedures are  analogous to the one $i=1$.

 First we write $A[G]=\coprod\limits_{g\in G} Ag$, as an $A^e$- module (or $A-A$ bimodule), and  so we get $HH^i(A[G])=\Ext_{A[G]^e}^i(A[G],A[G])\cong \coprod\limits_{(g,h)\in  G\times G}\Ext^i_{A^e}(Ag, Ah)$, for each  $i\geq 0$.

1. The action in case $i=0$

For each $g, h$ in $G$, let $f_{(g,h)}\in\Hom_{A^e}(Ag, Ah)$ and $(\sigma, \tau)\in G\times G$.  Then we \linebreak
consider the element
$(\sigma, \tau)(f_{(g,h)})\in\Hom_{A^e}(A\sigma g\tau^{-1}, A\sigma h\tau^{-1})$ defined by\linebreak
 $(\sigma, \tau)(f_{(g,h)})(a\sigma g\tau^{-1})= a\sigma f(g)\tau^{-1}$, with $a\in A$. 

Clearly it defines an action  on $\Hom_{A^e}(A[G], A[G])=\coprod\limits_{(g,h)\in G\times G}\Hom_{A^e}(Ag, Ah)$ and it is the action considered in Proposition \ref{gmod}.i).

2. The action in case $i=1$.

Let $\xi_{(g,h)}: 0\rightarrow Ah\rightarrow L \rightarrow Ag\rightarrow 0$ be a representative of an element in $\Ext^1_{A^e}(Ag, Ah)$ and $(\sigma, \tau)\in G\times G$. By  applying the exact functor $^{^{(\sigma^{\!\!-1}, \,\tau^{\!\!-1})}}()$ to this exact sequence  and using the isomorphism $^{^{(\sigma^{\!\!-1}, \,\tau^{\!\!-1})}} \!\!Ah\stackrel{\widetilde{}}\ra A\sigma h\tau^{-1}$ we obtain the  element, denoted by $(\sigma, \tau)(\xi_{(g,h)})\in \Ext^1_{A^e}(A\sigma g\tau^{-1}, A\sigma h\tau^{-1})$. It is easy to verify that it really defines an action on $\coprod\limits_{(g,h)\in G\times G}\Ext^1_{A^e}(Ag, Ah)$, and it is the action mentioned in Proposition \ref{gmod}.ii).

\begin{obser}\label{ide}

 We observe, in particular, that the subspace $\coprod\limits_{g\in G}\Ext^i_{A^e}(Ag, Ag)$ of $\Ext_{A^e}^i(A[G], A[G])$  under the action of $G\times G$ is taken to itself.  

On the other hand we also note that, for each $i\geq 0$ and $g\in G$,  $\Ext^i_{A^e}(A,A)$ is canonically isomorphic  to $\Ext^i_{A^e}(Ag, Ag)$ (see \ref{Misogm0}.iii). So, with this identification, we can consider an action of $G$ on $HH^i(A)=\Ext^i_{A^e}(A, A)$ such  as the one given  by the following: for each $g\in G$ and $\xi\in HH^i(A)$, $g\cdot\xi:= (g,g)(\xi)$. Consequently, we obtain that $G$ acts on $HH^{\bullet}(A)$.

\end{obser}

 Now we can show the main result of this section which gives a relation between
 the Hochschild cohomology algebras of $A$ and $A[G]$.

\begin{teo}\label{main}

 Let $A$ be a $k$-algebra and $G$ be a finite group acting on $A$. Then $G$ acts on the Hochschild cohomology $k$-algebra $HH^{\bullet}(A)$, and there is a ring monomorphism: $HH^{\bullet}(A)^G\hookrightarrow HH^{\bullet}(A[G]).$

\end{teo}

\begin{Demo}

First we write $A[G]=\coprod\limits_{g\in G}Ag$. Then  we  remark again that the action of $G\times G$ on $HH^i(A[G])$ ($i\geq 0$) enables having 
$HH^i(A[G])\cong (\Ext^i_{A^e}(A[G], A[G]))^{G\times G}\cong (\coprod\limits_{(g,h)\in G\times G}\Ext^i_{A^e}(Ag, Ah))^{^{G\times G}}$. So,
 it suggests us to identifying an element $\xi\in\Ext^i_{A^e}(A[G], A[G])$ with a matrix $\xi=(\xi_{(g,h)})_{g,h}$, with $\xi_{(g,h)}\in\Ext^i_{A^e}(Ag, Ah)$.

Now  we also remark that $HH^i(A)=\Ext^i_{A^e}(A,A)\cong\Ext^i_{A^e}(Ag,Ag)$, for any $g\in G$. So, according to Remark \ref{ide},  $G$  acts on $HH^i(A)$, and consequently on $HH^{\bullet}(A)$,  as it was indicated there. 

  Therefore the morphism we are looking for can be defined as: for each $i\geq 0$,  given $\xi\in HH^i(A)^G$ we take the element $\theta^i(\xi)\in  \coprod\limits_{(g,h)\in G\times G}\Ext_{A^e}^i(Ag, Ah)$ whose matrix representation $\theta^i=(\theta^i(\xi)_{(g,h)})_{g,h}$ is such that:
$$ \theta^i(\xi)_{(g,h)} = \cases{ 0 &if $g\neq h$\cr
 \xi &if $g = h$.\cr }$$
 Since $\xi\in (HH^i(A))^G$, then, by construction,  $\theta^i(\xi)$ is invariant under  the action of $G\times G$, and therefore the  map $\theta^i: HH^i(A)^G\to HH^i(A[G])$ is defined; and it is not difficult to verify that $\theta:HH^{\bullet}(A)^G\to HH^{\bullet}A[G]$, with $\theta =\oplus_i\theta^i$ is a monomorphism of rings. 
\qed  
\end{Demo}

\section {\sf Galois covering and Hochschild cohomology.}\label{covering}

In this section we are going to apply  the main theorem of the last section (Theorem \ref{main}) to show that also there is  a ring monomorphism 
from $HH^{\bullet}(\widetilde{A})^G$ into $HH^{\bullet}(A)$, where $\widetilde{A}$ is the  covering algebra of a $G$-graded $k$-algebra $A$.

 In \cite{CM} it was proven that, for a $k$-algebra $A$ graded by a finite group $G$,  the smash product   $A \# kG^{*}$ plays the role for graded rings that the skew group algebra $A[G]$ plays for group actions. So, in order to obtain the new pretended monomorphism we shall use the notion of smash product $A \# (kG)^{*}$, for  showing  the existence of an isomorphism between this product and the covering algebra of $A$ defined by $G$ and we apply the duality Theorem  3.5 in  \cite{CM}.

We recall here the definition  of covering algebra associated to a graded algebra. This definition was introduced in a preliminary version of \cite{GMS},  and it can be found in \cite{edu}. The definition of covering algebra  coincides with the one given by Green in \cite{G} for   quotient
 of path  algebras of quivers.

\begin{defi}\label{cov} Let $G$ be a finite group and $A=\coprod\limits_{g\in G}A_g$  be a $G$-graded $k$-algebra, with $A_g$ indicating the $k$-subspace of  the homogeneous elements of degree $g$.  The covering $k$-algebra associated to $A$, with respect to the given grading, denoted by
 $\widetilde A$, is  defined as follows. 
 As $k$- vector space \\
$\widetilde{A}= \coprod\limits_{(g,h)\in G\times G}\tilde{A} [g, h]$, where 
$\widetilde{A}[g,h]= A_{g^{-1}h}$, and 
the  multiplication is  defined
 in the following way: if $ \gamma \in \widetilde{A}[g, h]$ and $\gamma'\in \widetilde{A}[g', h']$. The product is in $ \widetilde{A}[g, h']$ and it is defined by:
$$\gamma\gamma'  = \cases{ 0 &if $g'\neq h$\cr
 \gamma\gamma' & if  $g' = h$.\cr }$$

\end{defi}

\begin{obser}\label{act} We observe that $G$ acts freely  on $\widetilde{A}$,  where the action of an element $\sigma\in G$ consists in to take   an element in $\widetilde{A}[g,h] $ to the same element, but now considered as an element in $\widetilde{A}[\sigma g, \sigma h]$. Moreover the canonical
vector space epimorphism $F:\widetilde{A}\to A$ which takes $\widetilde{A}[g, h]$ to $A_{g^{-1}h}$, is such that
 $F\sigma = F$, for all $\sigma \in G$, and the orbit space is $A$. So $\widetilde{A}$ is a Galois covering defined by $G$  in the sense of Gabriel and others,(\cite{ga, mape}).
\end{obser}

Now we review the definition of smash product  and some facts related to it (see \cite{CM}).

\begin{defi} Let $G$ be a finite group and  $A=\coprod\limits_{g\in G}A_g$ be a $G$-graded algebra. Let $k[G]^*$ be the dual algebra of $k[G]$,   and its  natural $k$-basis $\{p_g| g\in G\}$; that is , for any $g\in G$ and $x=\sum_{h\in G} a_h h\in k[G]$, $p_g(x)=a_g\in k$, and $p_g p_h=\delta_{g,h}p_h$, where $\delta_{g,h}$ is the Kronecker delta. The smash product, denoted by $A\# kG^*$,  is the vector space $A\otimes_k k[G]^*$ with the multiplication  given by 
$(a\#p_g)(b\#p_h) = a b_{g h^{-1}}\#p_h$ (here $a\# p_g$ denotes the element $a\otimes p_g$).
\end{defi}

The next proposition was first proved by Green, Marcos and Solberg in a preliminary version of  \cite{GMS}.

\begin{prop}\label{isom}
Let $G$ a finite group and $A=\coprod_{g\in G} A_g$ be a $G$-graded algebra. Then the smash product $A\# kG^*$ and the covering algebra $\widetilde{A}$ of $A$ are isomorphic algebras.
\end{prop}
\begin{Demo}
Any $a\in A$ can be written uniquely as $a = \sum_{h\in G} a_h$, where $a_h\in A_h$.
So we can define the following map $\Psi :A\# k[G]^*\longrightarrow \widetilde{A}$ by
$$\Psi((\sum_{h\in G}a_h)\# p_g) = \sum_{h\in G}a_h \in \coprod_{h\in G}\tilde{A}[g^{-1}h^{-1}, g^{-1}].$$ It is not hard to show that $\Psi$ is a bijective homomorphism of algebras.
\qed
\end{Demo}

In Remark \ref{act}, we have seen  the group of grading of $A$ acts on the covering algebra $\widetilde{A}$.  Now we also note that using the isomorphism in Proposition \ref{isom} we get a corresponding  action  on the smash product $A\#  kG^*$, which is given by $g(a\# p_h)=a\# p_{hg}$; and it coincides with the  one defined in Lemma 3.3 in \cite{CM}.

 With these remarks,  as a consequence of Theorem \ref{main}, of the isomorphism above and  the duality coactions of Cohen-Montgomery (Th.3.5 in \cite{CM}) we obtain the following proposition.

\begin{prop}\label{galois} Let $G$ be a finite group and $A$ be a $G$-graded $k$-algebra. Let  $\widetilde{A}$ be the  covering  algebra defined by the grading. Then $G$ acts on $HH^{\bullet}(\widetilde{A})$ and there is a ring monomorphism
from  $(HH^{\bullet}(\widetilde{A}))^G$ into  $HH^{\bullet}(A)$
\end{prop}
\begin{Demo}
As we have seen,  in Remark \ref{act}, $G$ acts on $\widetilde{A}$ as a group of automorphisms.  Then, on the one side,  from  Theorem \ref{main} it follows that  $G$ also acts on $HH^{\bullet}(\widetilde{A})$ and there is a monomorphism from from $HH^{\bullet}(\tilde{A})^G $ to $HH^{\bullet}((\widetilde{A})[G])$. But, on the other side, by Proposition \ref{isom}, $\tilde{A}$ and $A\# kG^*$ are isomorphic, and  according to our remark above this isomorphism leads  $G$ to act  on the smash product $A\# G$. So, by applying the duality theorem  for coactions (Th. 3.5 in \cite{CM}), we get  that $(A\# kG^*)[G]$ is isomorphic to the matrix ring $M_{|G|}(A)$ where $|G|$ 
denotes the order of the group $G$. Since Hochschild cohomology is an invariant by Morita equivalence
(in reality is  an invariant of derived equivalence, \cite{R, W}), then  it follows that $HH^{\bullet}(\widetilde{A}[G])$ is isomorphic to $HH^{\bullet}(A)$, and the proposition is  proved.
\qed
\end{Demo} 
\section {\sf The Hochschild cohomology of Koszul algebras}\label{Koszul}

In this section we discuss some facts related to  the Hochschild cohomology of Koszul algebras. In particular the ones concerning   to the preprojective algebras of Euclidean-type and to Auslander algebras of standard algebras, which are Koszul algebras. 

In order to study the Hochschild cohomology of Koszul algebras of finite global dimension we introduce   the construction of Koszul  complex for quadratic algebras. We are using a  similar procedure as was done  in \cite{Be} and \cite{gmmp} for generalized Koszul algebras (or $d$- Koszul algebras). According to our comments  in the introduction of this article, we use it for 2-Koszul algebra or simply   Koszul algebras. 
So we review some definitions and fix some notations.

 Let $k$ be a field and $Q$ be  a  finite quiver $Q$. We denote by  $kQ$  the path algebra of $Q$ and we indicate by $kQ_0$ the $k$-subalgebra whose underlying  vector space is the subspace generated by the vertex set $Q_0$ of $Q$. If  $Q_i$ is the set of paths of length $i$, then we denotes  $kQ_i$ the subspace of $kQ$ generated by $Q_i$. It is worth to note that this subspace is a $kQ_0$ bimodule. In this way, we will  consider the path algebra $kQ =\oplus_{i \geq 0} kQ_i$ as a graded algebra with the grading given by the length of the paths.

 Let $A = kQ/I$ where $Q$ is a finite quiver and $I$ is a two side ideal of $kQ$  generated by a set of quadratic relations (the $k$-algebra $A$ is called a quadratic algebra). Denoting by  $R$ the  set of homogeneous elements of degree two in I that  it is viewed as a $kQ_0$ sub-bimodule contained in $kQ_2$.

For each $n\geq 2$, we define  $K_n = \bigcap_{r +s + 2 =n} kQ_r.R.kQ_s$. Now we consider the following $A^e$-modules: 
$Q^i = 0$, if $ i< 0$;
$Q^0 = A\otimes_{_{kQ_0}} A$; $Q^1 = A\otimes_{_{kQ_0}} kQ_1\otimes_{_{kQ_0}}A$,
   and,  for  
$ n\geq 2$, $Q^n = A\otimes_{_{kQ_0}} K_n \otimes_{_{kQ_0}} A.$

It is clear that each $Q^i$ is an $A^e$-module finitely generated and projective. Observe also 
that each $Q^i$ is a submodule of $A\otimes_{_{kQ_0}} kQ_i\otimes_{_{kQ_0}} A$,
since $K_i$ is contained in $kQ_i$, for $i\geq 2$. 

 Now we construct, for each $n\geq 2$, the following $A^e$-morphisms \linebreak $f_n:A\otimes_{_{kQ_0}} kQ_n\otimes_{_{kQ_0}} A \to
A\otimes_{_{kQ_0}}kQ_{n-1}\otimes_{_{kQ_0}} A$ given by the formula

$$ f_n(a\otimes \alpha_1\cdots\alpha_n\otimes b) =
a \alpha_1\otimes \alpha_2\cdots \alpha_n\otimes b + (-1)^{n}
a\otimes\alpha_1\cdots\alpha_{n-1}\otimes \alpha_n b. $$

\begin{defi}
Let $A=kQ/I$  and $Q^i$ be the  $A^e$-modules  as above. Let  $d_i: Q^i\rightarrow Q^{i-1}$ be the maps such that $d_i =O$, for $i\leq 0$;  $d_1(1\otimes_{_{kQ_0}}\alpha\otimes_{_{kQ_0}} 1)= \alpha\otimes_{_{kQ_0}} 1 - 1\otimes_{_{kQ_0}} \alpha$, with $\alpha\in Q_1$,  and, for $n\geq 2$, $d_n$ is the restriction of $f_n$ to $Q^n$. It is very easy to see that  $d\circ d = 0$.  Then the complex $K^*(A)= ((Q^i)_i, (d_i)_i)$ is called the Koszul  complex of $A$.
\end{defi}

With the definition of  Koszul complex on the hands, we can utilize it for  characterizing the  Koszul algebras.  In order to get it we take the {\it augmented Koszul complex of $A$}:

$\cdots\ra Q^n\stackrel{d_n}{\ra}Q^{n-1}\stackrel{d_{n-1}}{\ra}\cdots\ra  
Q^1\stackrel{d_1}{\ra}Q^0\stackrel{d'_0}\ra A\ra 0$ where $d_n$, $n> 0$, is  as in $K^*(A)$ and $d'_0(a\otimes_{_{kQ_0}} b) = ab$  

\begin{teo}\cite{Be, gmmp}\label{Kos}
Let $A = kQ/I$ be a quadratic algebra.  The augmented  Koszul complex
$\cdots\ra Q^n\stackrel{d_n}{\ra}Q^{n-1}\stackrel{d_{n-1}}{\ra}\cdots\ra  
Q^1\stackrel{d_1}{\ra}Q^0\stackrel{d'_0}\ra A\ra 0$  is exact   if and only if $A$ is a Koszul algebra. 
\end{teo}

In case  $A=kQ/I$ is  a Koszul algebra, the augmented  Koszul complex of $A$ is a minimal graded projective resolution of $A$ as $A^e$-module. Furthermore, if $I$ is  an admissible ideal, then this resolution is also  a minimal projective resolution of $A$ in $A^e$-mod. 

The augmented  Koszul complex can be used to determine a lower bound for dimension of $HH^n(A)$, when $A$ is a Koszul algebra of global dimension $n$.\ It is obtained in the corollary below.

\begin{coro}\label{naoseanula} Let $A=KQ/I$ be a Koszul algebra of global dimension $n$. For each vertex $v\in Q_0$, $e_v$ denotes the associated idempotent  of $A$. Then 
$$\dim (HH^n(A)) \geq \dim (\coprod_{v\in Q_0}( e_vK_ne_v)) $$.
\end{coro}
\begin{Demo}
Since $A$ is a Koszul algebra and has the global dimension equal $n$, by Theorem   \ref{Kos}  we have, using the notations fixed above, that the long exact sequence 

\noindent $0\ra Q^n  \stackrel{d_n}\ra Q^{n-1}\cdots \ra Q^2  \stackrel{d_2}\ra Q^1\stackrel{d_1}\ra Q^0 \stackrel{d'_o}\ra A\ra 0$ is a graded projective resolution of $A$ in $A^e$-mod.  Then the Hochschild cohomology of $A$ can be computed as the cohomology groups of the complex:

\noindent $$0\rightarrow \Hom_{A^e}(A\otimes_{A_0}A, A)\stackrel{d_1^*}{\ra}\Hom_{A^e}(A\otimes_{A_0} kQ_1\otimes_{A_0} A, A)\ra\cdots\ra $$
\noindent $$ \Hom_{A^e}(A\otimes_{A_0} K_{n-1}
\otimes_{A_0} A, A)\stackrel{d_n^*}{\rightarrow} \Hom_{A^e} (A\otimes_{A_0}K_{n}\otimes_{A_0} A, A)\rightarrow 0 \cdots $$
\noindent where $A_0=kQ_0$.
 
On the other hand, it is easy to verify that  $\Hom_{A^e}(A\otimes_{A_0}A, A)\cong\Hom_{A_0^e}(A_0,A)\cong \coprod\limits_{v\in Q_0}e_vAe_v$, $\Hom_{A^e}(A\otimes_{A_0}kQ_1\otimes_{A_0} A, A)\cong\Hom_{A_0^e}(kQ_1,A)$, and, for $j\geq2$, $\Hom_{A^e}(A\otimes_{A_0}K_j\otimes_{A_0}A, A)\cong \Hom_{A_0^e}(K_j, A)$. Hence the last  complex is isomorphic to the following one:

$$ 0\rightarrow \coprod_{v\in Q_0} e_vAe_v\stackrel{d_1^*}{\ra}\Hom_{A_0^e}(kQ_1,A)\cdots\Hom_{A_0^e}(K_{n-1}, A)\stackrel{d_n^*}{\rightarrow} \Hom_{(A_0)^e}(K_n, A)\rightarrow 0, $$ 

\noindent where we also are denoting by  $d_i^*$ the induced maps by the isomorphisms mentioned above.

 We observe that the vector space  $\Hom_{(A_0)^e}(K_j, A)$, for $j\geq 2$,  can be naturally graded by the induced grading of $A$; that is, we say  that a map $f$ is homogeneous of degree $t$ if its image is contained in homogeneous component of degree $t$ of $A$. So, it easy to see that the last complex  is a complex of graded vector spaces and that  the image of each $d_j^*$ is contained in the direct sum of the homogeneous subspaces of degree bigger than zero. Then, for $j\geq 2$, the image of $d_j^*$ does not intersect the degree zero   component $(\Hom_{A^e_0}(K_j, A))_0\cong \coprod_{v\in Q_0}(e_vK_je_v)$ of $\Hom_{A^e_0}(K_j, A)$. In particular,  the degree zero component   
$(\Hom_{A^e_0}(K_n, A))_0\cong\coprod_{v\in Q_0}(e_vK_ne_v)$ does not intersect the image of $d_n^*$  and since $HH^n(A)=\Coker d^*_n$, a simple computation of dimensions shows our statement.
 \qed
\end{Demo}
 
 Among the Koszul algebras we are going to point out the $\mathbb{C}$-preprojective algebras of Euclidean type and the Auslander algebra of a standard, representation finite-type $k$- algebra. We will see that these algebras are Koszul algebras and as consequence of it, via Corollary \ref{naoseanula}, we get interesting datum  for $HH^2(A)$.

First let us review the definition of  preprojective algebras.

\begin{defi} Let $Q$ be a finite quiver and $k$ a field. 
 Consider the quiver $\hat Q$ whose vertex set $\hat Q_0= Q_0 $ and the arrows set $\hat{Q_1}= Q_1\cup Q_1^{op}$, where $Q^{op}$ denotes the opposite quiver of $Q$.  For each arrow $\alpha\in Q_1$ we write $\hat{\alpha}$ for the corresponding arrow in the opposite quiver. 
The preprojective $k$-algebra associated to $Q$ (or briefly the preprojective $k$-algebra of $Q$), denoted by ${\mathcal P}(Q)$, is 
is the $k$-algebra $k\hat{Q}/I$, where $I$ is the ideal generated  by the relations $\sum\limits_{\alpha\in Q_1}  \alpha\hat\alpha$  and  $\sum\limits_{\alpha\in Q_1} \alpha\hat\alpha$.
\end{defi}

We remark that it is well-known that the preprojective algebra constructed as above only depends on the underlying graph of the quiver $Q$; that is, quivers having the same underlying graph define isomorphic preprojective algebras.

The Hochschild cohomology of preprojective algebras associated to  Dynkin diagrams $A_n$ were studied in \cite{ES}. 

We mention the following result about preprojective $\mathbb{C}$-algebras associated to  Euclidean diagrams (see \cite{Ch, Le, ReB}).

\begin{teo}\label{len}
The preprojective $\mathbb{C}$-algebras associated to an  Euclidean diagram  are Morita equivalent to the skew group
algebras  $\mathbb{C}[x,y][G]$, with $G$ a polyhedral group.
\end{teo}

So this theorem can be used in order to study some properties of the preprojective $\mathbb{C}$-algebra of 
Euclidean-type through properties of certain skew group rings. For instance, from this theorem we obtain  that a preprojective $\mathbb{C}$-algebra associated to   Euclidean diagrams have global dimension 2 (recall that  $gldim(\mathbb{C}[x,y][G])= gldim(\mathbb{C}[x,y]) = 2$). Moreover, since the preprojective algebras are always quadratic algebras, then we also get that preprojective $\mathbb{C}$-algebras of Euclidean-type are Koszul algebras.

In this  point of view we obtain as a consequence of  Corollary \ref{naoseanula}, the following fact about the second Hochschild cohomology group of  preprojective $\mathbb{C}$-algebras of Euclidean-type. 

\begin{coro}
 Let $A$ be a preprojective $\mathbb{C}$-algebra associated to an  Euclidean diagram. Then $HH^2(A) \neq 0$
\end{coro}
\begin{Demo} We have that $A=\mathbb{C}{\hat Q}/I$ where $Q$ is an Euclidean diagram. As we have commented above $A$ is a Koszul algebra of global dimension two. Since $e_vK_2e_v$ is not zero,  for any vertex $v$, the statement follows from the corollary \ref{naoseanula}.
\qed
\end{Demo}

We note that Theorem \ref{len} can be used this  to describe the centre of preprojective $\mathbb{C}$-algebras of Euclidean-type, once their centres are  the ones of the skew group rings $\mathbb{C}[x,y][G]$, for suitable groups $G$. Then it seems interesting to study  the Hochschild cohomology of the skew group ring, in particular its centre.  In order for studying it the following lemma will be useful.

\begin{lem}\label{leM} 
Let $R$ be a commutative ring and $G$ be a group acting  on $R$. Let $g $ be an element in $G$ and suppose that  there is 
$\alpha \in R$ such that $\alpha -g(\alpha) $  is not a zero divisor in $R$. Then
$\Hom_{R^e}(R, Rg) = 0$.
\end{lem}
\begin{Demo}
Let $f\in\Hom_{R^e}(R, Rg)$. Then $\alpha f(1) = f(1).\alpha = f(1)g(\alpha) = g(\alpha)f(1)$, and it implies that $(\alpha - g(\alpha))f(1) = 0$. Since $\alpha - g(\alpha)$ is not a zero divisor in $R$, it follows that $f(1) =0$, and consequently $f =0$.
\qed
\end{Demo}
 
Now we be able to describe the center of a skew group ring.

\begin{prop}\label{cen}Let $R$ be a commutative domain and $G$ be a finite group of automorphism of  $R$. Then  
$\centro(R[G])$ is isomorphic to $R^G$.
\end{prop}

\begin{Demo}

First we observe that $Hom_{R^e}(Rg, Rh)\cong \Hom_{R^e}(R, Rhg^{-1})$, for any $g, h$ in $G$. So, since $G$ acts faithfully on $R$ and $R$ is a domain,  by Lemma \ref{leM} we obtain that $\Hom_{R^e}(Rg, Rh)=0$, for any 
  $g\neq h$ in $G$.

Now, if we write  $R[G]= \coprod\limits_{g\in G} Rg$ as $R^e$-module, 
 then  we have that $\centro (R[G])=HH^0(R[G])=\Hom_{R[G]^e}(R[G], R[G]) = (\Hom_{R^e}(R[G], R[G]))^{G\times G} =\linebreak
(\coprod\limits_{(g,h)\in G\times G}Hom_{R^e} (Rg, Rh))^{G\times G} \cong (\coprod\limits_{g\in G} Hom_{R^e}(Rg, Rg))^{G\times G}$. 

Recalling the action defined in section 2 and Remark \ref{ide} we have that
\linebreak
 $ (\coprod\limits_{g\in G} Hom_{R^e}(Rg, Rg))^{G\times G}\cong (Hom_{R^e}(R,R))^G\cong R^G$, and the statement is proved.

\qed
\end{Demo}

We remark that the proof of Proposition \ref{cen}   could be obtained by a direct computation, but we have optioned by the proof above  for illustrating how to use our methods.  

For the next corollary we need some additional terminology. We are going to consider the Auslander algebra associated to a $k$-algebra of representation finite type. So, we recall
 the definition of  Auslander algebras.

Let $A$ be a $k$-algebra  of representation finite type (i.e. up to isomorphism there exist only finitely many indecomposable $A$-modules). Let $X_1, X_2, \cdots X_m$ be a list of representatives from the isomorphism classes of indecomposable $A$-modules and let $X=\oplus_i X_i$. The $k$-algebra $\Lambda=\mbox{End}_A (X)$ is called {\it Auslander algebra of $A$}. Recall that $A$ is  said to be {\it standard} if $\Lambda$ is isomorphic to the quotient of
 the path algebra $k\Gamma_A$ of the Auslander-Reiten quiver $\Gamma_A$ of $A$ by the ideal generated by the mesh relations. 

 We denote by mod-$A$ the category of finitely generated left $A$-modules, by $\ind  A$ the subcategory of mod-$A$ with one representative of each isoclass of indecomposable $A$-module and by $\tau_A$ the Auslander-Reiten translate DTr.

\begin{coro}
Let $A$ be a standard representation-finite type  $k$-algebra and $\Lambda$ be its   Auslander algebra. Then 
$\dim HH^2(\Lambda) \geq \#\{M\in \ind A: \tau_A M = M\} $
\end{coro}
\begin{Demo}
It is known that the Auslander algebra $\Lambda$ of a representation-finite algebra $A$   has global dimension two. Moreover, since $A$ is standard, we have that $\Lambda\cong k\Gamma_A/I$, where $\Gamma_A$ is the Auslander-Reiten quiver  of $A$ and $I$ is the ideal generated by the mesh relations (so quadratic relations). Hence $\Gamma$ is a Koszul algebra.

By construction of the quiver $\Gamma_A$ and by the conditions on $I$, it is clear that the number of elements of the set $\{M\in \ind A:\tau_A M) = M\} $ is the dimension of the degree zero component of $Hom_{k(\Gamma_A)_0^e}(K_2, \Lambda)$. But that component is isomorphic to  $\coprod_{v\in (\Gamma_A)_0}e_vK_2e_v$, and therefore the result follows from corollary \ref{naoseanula}
\qed
\end{Demo}

\medskip

\noindent Departamento de Matem\'atica - IME, Universidade de S\~ao Paulo,\\
 C. Postal 66281, CEP 05315- 970, S\~ao Paulo, SP, Brasil\\
{email: enmarcos@ime.usp.br}

\medskip
\noindent Instituto de Matematicas, UNAM- Campus Morelia,\\
 Apartado Postal 61-3, CP 58089, Morelia, Michoac\'an, Mexico\\
{email: mvilla@matmor.unam.mx}

\medskip
\noindent Departamento de Matem\'atica - IME, Universidade de S\~ao Paulo,\\
 C. Postal 66281, CEP 05315- 970, S\~ao Paulo, SP, Brasil\\
{email: bel@ime.usp.br}

\end{document}